\documentclass[onecolumn, notitlepage, superscriptaddress]{revtex4-1}

\usepackage{graphicx}
\usepackage[T1]{fontenc}    
\usepackage[bookmarks=false,hidelinks]{hyperref} 
\usepackage{url}            
\usepackage{amsfonts}       
\usepackage{amsmath,amssymb}
\usepackage{amsthm}
\usepackage{bbm}
\usepackage{physics}
\usepackage{mathtools}
\usepackage{latexsym}
\usepackage{subcaption}
\usepackage{enumitem}
\usepackage[final]{changes}
\usepackage{xcolor}
\usepackage{here}
\usepackage[ruled]{algorithm2e}

\usepackage{academicons}
\newcommand{\orcid}[1]{\href{https://orcid.org/#1}{\textcolor[HTML]{A6CE39}{\aiOrcid}}}

\theoremstyle{definition}
\newtheorem{theorem}{Theorem}
\newtheorem{corollary}{Corollary}
\newtheorem{proposition}{Proposition}
\newtheorem{definition}{Definition}

\newtheorem{assumption}{Assumption}

\newtheorem{remark}{Remark}

\usepackage{cleveref}

\crefname{equation}{}{}
\crefname{figure}{Fig.}{Figs}
\crefname{table}{Table}{Tables}
\crefname{algocf}{Algorithm}{Algorithms}
\Crefname{algocf}{Algorithm}{Algorithms}
\crefname{algorithm}{Algorithm}{Algorithm}
\crefname{section}{Section}{Sections}

\crefname{theorem}{Theorem}{Theorems}
\crefname{lemma}{Lemma}{Lemmas}
\crefname{corollary}{Corollary}{Corollaries}
\crefname{assumption}{Assumption}{Assumptions}
\crefname{definition}{Definition}{Definitions}
\crefname{remark}{Remark}{Remarks}
\crefname{proposition}{Proposition}{Propositions}
\crefname{example}{Example}{Examples}

\DeclareFontFamily{U}{mathx}{}
\DeclareFontShape{U}{mathx}{m}{n}{<-> mathx10}{}
\DeclareSymbolFont{mathx}{U}{mathx}{m}{n}
\DeclareMathAccent{\widehat}{0}{mathx}{"70}
\DeclareMathAccent{\widecheck}{0}{mathx}{"71}

\DeclarePairedDelimiterX{\inp}[2]{\langle}{\rangle}{#1, #2}

\DeclareMathOperator*{\argmin}{argmin}

\newcommand{\calH}{{\mathcal H}}

\newcommand{\calN}{{\mathcal N}}

\newcommand{\calX}{{\mathcal X}}

\newcommand{\bbE}{{\mathbb E}}

\newcommand{\bbN}{{\mathbb N}}

\newcommand{\bbR}{{\mathbb R}}

\newcommand{\bbV}{{\mathbb V}}

\newcommand{\rmd}{{\mathrm d}}

\newcommand{\tilk}{{\widetilde k}}

\newcommand{\tilu}{{\widetilde u}}
\newcommand{\tilv}{{\widetilde v}}

\newcommand{\tilN}{{\widetilde N}}

\begin{document}

\title{
An Uncertainty-aware, Mesh-free Numerical Method for Kolmogorov PDEs
}

\author{Daisuke Inoue}
\email{daisuke-inoue@mosk.tytlabs.co.jp}
\affiliation{%
Toyota Central R\&D Labs., Inc.\\
Nagakute, Aichi 480-1192, Japan
}%
\author{Yuji Ito}
\affiliation{%
Toyota Central R\&D Labs., Inc.\\
Nagakute, Aichi 480-1192, Japan
}%
\author{Takahito Kashiwabara}
\affiliation{%
Graduate School of Mathematical Sciences, the University of Tokyo\\
3-8-1 Komaba, Meguro-ku, Tokyo 153-8914, Japan
}%
\author{Norikazu Saito}
\affiliation{%
Graduate School of Mathematical Sciences, the University of Tokyo\\
3-8-1 Komaba, Meguro-ku, Tokyo 153-8914, Japan
}%
\author{Hiroaki Yoshida}
\affiliation{%
Toyota Central R\&D Labs., Inc.\\
Nagakute, Aichi 480-1192, Japan
}%

\begin{abstract}
  This study introduces an uncertainty-aware, mesh-free numerical method for solving Kolmogorov PDEs.
  In the proposed method, we use Gaussian process regression (GPR) to smoothly interpolate pointwise solutions that are obtained by Monte Carlo methods based on the Feynman--Kac formula.
  The proposed method has two main advantages: 1. uncertainty assessment, which is facilitated by the probabilistic nature of GPR, and 2. mesh-free computation, which allows efficient handling of high-dimensional PDEs. 
  The quality of the solution is improved by adjusting the kernel function and incorporating noise information from the Monte Carlo samples into the GPR noise model.
  The performance of the method is rigorously analyzed based on a theoretical lower bound on the posterior variance, which serves as a measure of the error between the numerical and true solutions. 
  Extensive tests on three representative PDEs demonstrate the high accuracy and robustness of the method compared to existing methods.
\end{abstract}

\maketitle

\section{Introduction}\label{sec:intro}

Kolmogorov partial differential equations (PDEs) have found wide application in physics~\cite{Pascucci2005Kolmogorov}, chemistry~\cite{van1992stochastic}, finance~\cite{black1973pricing}, and the natural sciences~\cite{edelstein2005mathematical}.
Various numerical methods such as finite difference~\cite{Brennan1978Finite,Kumar2006Solution,Kushner1976Finite,Zhao2007Compact} and finite element method~\cite{Floris2013Numeric,Zienkiewicz2013Finite} have been proposed and widely adopted by practitioners.
While these numerical methods are undoubtedly of high value in solving various engineering problems, two challenges remain in approaching more complex real-world problems: 

\begin{enumerate}
  \item
  Evaluating the validity of numerical results is challenging.
  Error analysis of existing numerical methods requires careful mathematical analysis, which is different for each PDE of interest.
  Furthermore, the tightness of these bounds is sometimes not enough to assess the detailed error of the computed results.
  \item 
  The computational complexity dramatically increases with the degrees of freedom in meshes. 
  Solving the linear equations arising from discretizing the PDE with meshes typically requires a computational complexity that scales cubically with the mesh size. 
  This growth, scaled to the powers of the spatial dimensions, imposes impractical time requirements for solving high-dimensional PDEs.
\end{enumerate}

The first challenge has long been the focus of error analysis in conventional numerical methods, such as finite difference and finite element methods~\cite{Strikwerda2004Finitea,Verfurth2013Posteriori}.
In particular, in the finite element method, a priori or a posteriori error estimation plays a crucial role, which is performed before and after the numerical solution is obtained, respectively~\cite{Chamoin2023Introductory,Ladeveze1983Error}.
A priori estimation primarily involves truncation error analysis~\cite{Chazan1968Introduction}, while a posteriori estimation often involves flux recovery methods, residual methods, and duality-based approaches~\cite{Strouboulis1992Recent}. 
In general, tight error analysis requires careful mathematical evaluation, which depends on the form of the PDEs.
It is therefore desirable to evaluate the error numerically, regardless of the form of the PDEs.
Furthermore, the reliance on meshes in these conventional methods poses computational challenges for solving high-dimensional equations.

To address the second issue, Monte Carlo-based mesh-free calculations are actively investigated.
In particular, methods based on the Feynman--Kac (FK) formula have attracted attention as a promising strategy for solving high-dimensional problems~\cite{Exarchos2016Learning,Exarchos2018Stochastic,Han2018Solving,Hawkins2021Time,LeCavil2019Forward,Pereira2020FeynmanKac,Pham2015FeynmanKac}.
This method selects a point in space-time to be solved, calculates multiple stochastic differential equations (SDEs) starting from the point, and approximates the solution of the PDE at the selected point as an ensemble average of the SDE solutions.
Here, the FK formula provides a solution at a single point in space-time, requiring regressions from multiple sampled points to derive a solution over the entire domain.
Regression via linear interpolation and neural networks~\cite{Beck2021Solving,Exarchos2018Stochastic,Richter2022Robust,Takahashi2023Solving} has been proposed.
Nonetheless, these approaches give only estimates of the PDE solutions and cannot guarantee the accuracy of the estimated results.

This study proposes a new numerical method to solve the above two challenges simultaneously.
In the proposed method, Gaussian Process Regression (GPR) is used to obtain solutions to PDEs over an entire domain while evaluating their uncertainty. 
Specifically, GPR smoothly interpolates solutions at discrete points obtained by a Monte Carlo method based on the FK formula.
GPR not only allows us to obtain regression results as a mean function of the posterior distribution but also enables us to assess the uncertainty of the regression through the posterior variance.
The latter allows a posteriori probabilistic evaluation of the error between the numerical solution and the true solution under appropriate assumptions.

GPR is compatible with Monte Carlo methods from two perspectives that are unattainable with other regression techniques.
Firstly, GPR can incorporate a priori knowledge about the regularity of PDE solutions.
In particular, designing a kernel function allows for determining the differentiability of the regression results.
Secondly, using the variance information of noise contained in the Monte Carlo samples increases the accuracy of the regression. 
The variance is easily estimated by the sample variance of the Monte Carlo samples.

The use of GPR allows not only a posteriori evaluation of the uncertainty of the numerical solution through posterior variance, but also a priori estimation of the posterior variance.
This is achieved through the following two theoretical analyses.
First, we establish a lower bound on the posterior variance which is crucial for estimating the best-case performance in numerical solutions. 
This bound is obtained by extending the asymptotic theory of standard GPR to the heteroscedastic case, where the noise term has different magnitudes at each point in space-time (see Theorem \ref{thm:convergence-MSE-hetero-min}).
Yet, this lower bound includes the variance of the noise term in GPR, making it difficult to evaluate before computing the posterior distribution. 
To address this, we probabilistically assess the error between the true variance and the sample variance, which can be easily calculated from the FK sample.
This is given as a concentration inequality using Chebyshev's inequality (see Theorem \ref{thm:err-bound}).
The combination of these two analytical results allows us to obtain a priori error bounds for the proposed numerical methods.

The remainder of the paper is structured as follows: In \cref{sec:fkgp-problem}, we define the class of the PDEs to be addressed.
In \cref{sec:fkgp-algorithm}, we propose a numerical method that combines the GPR and the FK formula.
Theoretical results for the estimation of the posterior variance in the proposed method are presented in \cref{sec:theory}.
In \cref{sec:numerics}, we apply the proposed method to three types of high-dimensional PDEs and show that more accurate and less uncertain solutions are obtained than with existing methods.
\cref{sec:fkgp-conclusion} provides a brief summary and suggests possible future work.

\section{Problem Formulation}\label{sec:fkgp-problem}

Throughout this paper, we consider the following Cauchy problem for the Kolmogorov PDEs:
\begin{align}
  &\begin{aligned}\label{eq:BPDE}
    \partial_t v(t,x)&+\frac{1}{2} \operatorname{tr}\left\{a(t, x) a(t, x)^{\top} \partial_{xx} v(t,x)\right\} + b(t, x)^\top \partial_x v(t, x)\\
    &\quad +c(t, x) v(t,x) + h(t, x)=0,\quad\quad\quad (t, x) \in[0, T) \times \mathbb{R}^d,
  \end{aligned}  \\
  &v(T,x)=g(x), \qquad\qquad\qquad\qquad\qquad\qquad\qquad\qquad\qquad x \in\mathbb{R}^d,\label{eq:BPDE-terminal}
\end{align}
where $T>0$ denotes the terminal time, $t\in[0,T]$ denotes the time, and $v:[0,T]\times \bbR^d\to\bbR$ denotes the variables of the PDE.
The function $b:[0, T] \times \mathbb{R}^d \rightarrow \mathbb{R}^d$ denotes the advection coefficient, $a:[0, T] \times \mathbb{R}^d \rightarrow \mathbb{R}^{d \times m}$ denotes the diffusion coefficient, $c:[0, T] \times \mathbb{R}^d \rightarrow \mathbb{R}$ denotes the reaction coefficient, $h: [0, T] \times \mathbb{R}^d \rightarrow \mathbb{R}$ denotes the generative coefficient, and $g:\mathbb{R}^d \rightarrow \mathbb{R}$ denotes the terminal condition.
We assume that each function satisfies the following conditions.

\begin{assumption}\label{Asmp:pde} 
  The functions $b:[0, T] \times \mathbb{R}^d \rightarrow \mathbb{R}^d, a:[0, T] \times \mathbb{R}^d \rightarrow \mathbb{R}^{d \times m}$ $c:[0, T] \times \mathbb{R}^d \rightarrow \mathbb{R}$, $h: [0, T] \times \mathbb{R}^d \rightarrow \mathbb{R}$ are uniformly continuous.
  The function $c$ is bounded.
  There exists a constant $L>0$ and the following are satisfied for $\varphi \in\{ b, a, c, h\}$:
  \begin{align}
    \begin{aligned}
      |\varphi(t, x)-\varphi(t, x^{\prime})| &\leq L|x-x^{\prime}|, \quad \forall t \in[0, T], x, x^{\prime} \in \mathbb{R}^d, \\
      |\varphi(t, 0)| &\leq L, \qquad\qquad \forall t \in[0, T].
    \end{aligned}
  \end{align}
\end{assumption}

\begin{proposition}[Ref.~\cite{Yong1999Stochastic}]
  Suppose that \cref{Asmp:pde} holds.
  Then, \cref{eq:BPDE,eq:BPDE-terminal} have a unique viscosity solution.
\end{proposition}

Finding a solution to \cref{eq:BPDE,eq:BPDE-terminal} in closed form is possible only in very limited cases and thus numerical calculations are usually required.
The goal of this study is to find a numerical solution $v$ of \cref{eq:BPDE,eq:BPDE-terminal}.
More specifically, under the known functions $a, b, c, h$, and $g$, we aim to find the solution $v(0,x)\ (x\in\calX)$ at the initial time $t=0$, where the set $\calX\subseteq\bbR^d$ is the domain for which we want to know the solution.
In addition, we also aim to assess the uncertainty of the solution, that is, how much the results of the numerical calculation deviate from the true solution.

\begin{remark}\label{rem:initial-terminal}
  The terminal value problem \cref{eq:BPDE,eq:BPDE-terminal} is transformed into the initial value problem by introducing the variable transformation $w(t,x)\coloneqq v(T-t,x)$: 
  \begin{align}
    &\begin{aligned}\label{eq:BPDE2}
      -\partial_s w(s,x)&+\frac{1}{2} \operatorname{tr}\left\{a(s, x) a(s, x)^{\top} \partial_{xx} w(s, x)\right\} + b(s, x)^\top \partial_x w(s,x)\\
      &\quad +c(s, x) w(s, x) + h(s, x)=0, \quad\quad\ (s, x) \in[0, T) \times \mathbb{R}^d,
    \end{aligned} \\
    &w(0,x)=g(x), \qquad\qquad\qquad\qquad\qquad\qquad\qquad\qquad\qquad\qquad x \in\mathbb{R}^d,\label{eq:BPDE2-initial}
  \end{align}
  Thus, for initial value problems, this variable transformation should be performed before applying the proposed method described below.
\end{remark}

\section{Proposed Algorithm}\label{sec:fkgp-algorithm}

In this section, we propose a method to find the solution $v(0,x)\ (x\in\calX)$ of the PDE \cref{eq:BPDE,eq:BPDE-terminal} at the initial time $t=0$.
We use GPR to compute the possibility that each candidate function is a solution $v$ of the PDE.
Using prior knowledge and data sets, the GPR provides such a possibility as a distribution over possible functions.
The prior knowledge represents the possibility regarding the solution $v$ as a probability distribution 
$p(v)$.
If the data set $U$ associated with $v$ is observed, the posterior probability 
$p(v|U)$ is analytically obtained by Bayesian estimation.
In this sense, the proposed data-driven method estimates the possibility $p(v|U)$
that each candidate of $v$ is the solution to the PDE instead of solving the PDE directly.
In \cref{sec:GPR}, we describe how to obtain the posterior distribution by using GPR.
In \cref{sec:FK}, we describe how to obtain the data set that approximates the solution of the PDE by using the Feynman--Kac formula.

\subsection{Gaussian Process Regression}\label{sec:GPR}

We design the prior as a Gaussian process.
For any $X=[x_1,\ldots,x_N]^\top\in\calX^N$, we assume that the probability distribution $p(v_X|X)$, where 
\begin{align}
  v_{X}\coloneqq[v(0,x_1),\ldots,v(0,x_N)]^\top\in\bbR^N,
\end{align}
is a normal distribution:
\begin{align}\label{eq:f_is_gp}
 p(v_X|X) = \calN(0, k_{XX}).
\end{align}
Here, the matrix $k_{XX}\in\bbR^{N\times N}$ is defined as 
$\left[k_{XX}\right]_{i j}=k\left(x_i, x_j\right)$ with a positive-definite kernel function $k:\calX\times\calX\to \bbR$, where $[k_{XX}]_{ij}$ denotes the $(i,j)$ component of $k_{XX}$.
See Appendix for a detailed definition of the Gaussian process.
The method for designing the kernel function $k$ is discussed in \cref{sec:FK}.

We assume that the data set is given as follows.
For $N$ observation points $X=\left[x_1, \ldots, x_{N}\right]^\top\in \mathcal{X}^N$ and $M$ sample size, we assume that the data
\begin{align}
  U=\left[u_1(x_1),\ldots,u_M(x_1), \ldots, u_1(x_{N}),\ldots, u_M(x_{N})\right]^\top \in \mathbb{R}^{MN},
\end{align}
which approximates the solution of the PDE, are obtained.
We assume that these data satisfy
\begin{align}\label{eq:regression_model}
  u_j(x_i)= v(0,x_i) + \xi_{i,j}, \quad i\in\{1, \ldots, N\},\ j\in\{1, \ldots, M\},
\end{align}
where $\xi_{i,j}\ (i\in\{1, \ldots, N\},\ j\in\{1, \ldots, M\})$ are random variables that are independent and follow a normal distribution with mean $0$ and variance $r(x_i)$:
\begin{align}\label{eq:noise-var-hetero}
  \xi_{i,j} \sim \mathcal{N}\left(0, r(x_i)\right), \quad i\in\{1, \ldots, N\},\ j\in\{1, \ldots, M\}.
\end{align}
The function $r:\calX\to \bbR_{\ge 0}$ is the noise variance function.
To reduce the effect of noise, we derive the posterior distribution using the sample mean $\overline u(x_i)$ for the data at the same observation point $x_i$:
\begin{align}\label{eq:mean_data_set}
  \overline U = \left[\overline u(x_1),\ldots,\overline u(x_N)\right]^\top\in\bbR^N,
\end{align}
where $\overline {{u}} (x_i)\coloneqq \sum_{j=1}^M {u}_j (x_i)/M$.
Since the variance of the sample mean is the variance of each sample divided by $M$, \cref{eq:regression_model} is rewritten as follows:
\begin{align}\label{eq:mean_regression_model}
  \overline u(x_i) = v(0,x_i) + \overline \xi_{i}, \quad i\in\{1, \ldots, N\},
\end{align}
where
\begin{align}\label{eq:mean_noise-var-hetero}
  \overline \xi_{i}\sim \mathcal{N}\left(0, \frac{r(x_i)}{M}\right), \quad i\in\{1, \ldots, N\}.
\end{align}
In summary, the dataset is observed to satisfy
\begin{align}\label{eq:observed_data_set}
  p(\overline U|v_X,X) = \calN\left(v_X, \frac{r_{XX}}{M}\right),
\end{align}
where $r_{X X}\in \mathbb{R}^{N \times N}$ is a noise variance matrix of the form $\left[r_{X X}\right]_{i j}=\delta_{ij} r(x_i)$.
The methods for obtaining the data set $\overline U$ and the noise variance matrix $r_{XX}$ are described in \cref{sec:FK}. 

In GPR, the posterior distribution can be analytically obtained from the prior and the data set, as stated in the next proposition.
Recalling that $k$ is the kernel function representing the prior in \cref{eq:f_is_gp}, we define $k_{X x}=k_{x X}^{\top}=\left[k\left(x_1, x\right), \ldots, k\left(x_{N}, x\right)\right]^{\top}$.

\begin{proposition}[Ref.~\cite{Goldberg1997Regression}]\label{thm:posterior-hetero}
  Suppose that \cref{eq:f_is_gp,eq:observed_data_set} holds.
  For any $x\in\calX$, 
  the posterior distribution 
  is calculated as
  \begin{align}\label{eq:posterior_distribution}
    p(v(0,x)| X,\overline U,x) = \calN(\tilu(x), \tilk(x,x)),
  \end{align}
  where
  \begin{align}
    \tilu(x)&=k_{x X}\left(k_{X X}+\frac{r_{X X}}{M}\right)^{-1}\overline U, \quad x \in \mathcal{X},\label{eq:GPR_mean-hetero}\\
    \tilk\left(x, x^{\prime}\right)&=k\left(x, x^{\prime}\right)-k_{x X}\left(k_{X X}+\frac{r_{XX}}{M}\right)^{-1} k_{X x^{\prime}}, \quad x, x^{\prime} \in \mathcal{X}. \label{eq:posterior_variance}  
  \end{align}
\end{proposition}

Here, we estimate the solution of the PDE by evaluating the posterior mean function \cref{eq:GPR_mean-hetero}.
We use the posterior variance $\widetilde\sigma^2(x) \coloneqq \tilk(x,x)$ to assess the uncertainty of the estimates.
This value represents the expectation of the mean squared error (MSE) between the true solution $v(0,x)$ and the posterior mean function, where the expectation is evaluated for the probability $p(v(0,x)| X,\overline U,x)$.

\begin{remark}
  In standard GPR, the random variable $\xi_{i,j}$ has constant variance $r(x)\equiv \sigma^2$, which is a hyperparameter learned to maximize the likelihood.
  As described in \cref{sec:FK}, in our method the noise variance is estimated as a function over the space $\calX$ by using the sampled data set.
  Although the concept of heteroscedastic GPR has already been studied for many years~\cite{lazaro2011variational,le2005heteroscedastic,Liu2020Largescale,Makarova2021Riskaverse}, to the best of our knowledge, this is the first time that it has been used for the numerical calculation of PDEs.
\end{remark}

\subsection{Prior and Data Set}\label{sec:FK}

In the GPR introduced in the previous section, the following design parameters have not yet been determined:
\begin{enumerate}
  \item The kernel function $k$ in \cref{eq:f_is_gp}
  \item The data set $U$ and the noise variance matrix $r_{XX}$ that satisfy \cref{eq:observed_data_set}
\end{enumerate}

We first discuss the kernel function $k$.
We propose the use of the following Mat\'{e}rn kernel to reflect a priori knowledge about the regularity of the PDE solution:
\begin{align}\label{eq:Matern}
  k\left(x, x^{\prime}\right)=\frac{1}{2^{\alpha-1} \Gamma(\alpha)}\left(\frac{\sqrt{2 \alpha}\left\|x-x^{\prime}\right\|}{h}\right)^\alpha K_\alpha\left(\frac{\sqrt{2 \alpha}\left\|x-x^{\prime}\right\|}{h}\right),
\end{align}
where $\alpha>0$ and $h>0$ are hyperparameters, $\Gamma$ is the Gamma function, and $K_\alpha$ is the modified Bessel function of the second kind of order $\alpha$.

This is justified by the following.
According to the Moore--Aronszajn theorem~\cite{Aronszajn1950Theory}, for any positive-definite kernel $k$, there is only one reproducing kernel Hilbert space (RKHS) $\mathcal{H}_k$ for $k$.
The norm defined on this RKHS can evaluate the smoothness of the function.
Particularly, the RKHS of the Mat\'{e}rn kernel is known to be norm-equivalent to the following  Sobolev space:
\begin{proposition}[RKHSs of Mat\'{e}rn kernels: Sobolev spaces~\cite{Kanagawa2018Gaussian,Wendland2004Scattereda}]
  Let $k$ be the Mat\'{e}rn kernel on $\mathcal{X} \subseteq \mathbb{R}^d$.
  Suppose that $s\coloneqq\alpha+d / 2$ is an integer.
  Then, the RKHS $\mathcal{H}_{k}$ is norm-equivalent to the
  Sobolev space of order $s$.
\end{proposition}

Therefore, in situations where information about the regularity of the PDE solutions is known, the smoothness of the regression results can be made to match the smoothness of the true solution by appropriately choosing the parameter $\alpha$ in the Mat\'{e}rn kernel.

Next, to observe the data that satisfy \cref{eq:regression_model}, we employ the Feynman-Kac formula.
Using the coefficients of \cref{eq:BPDE,eq:BPDE-terminal}, we define a random variable
\begin{align}
  \begin{aligned}\label{eq:FK_one_sample}
    \widehat u(t,x)&\coloneqq \int_t^T h(s, X(s ; t, x)) e^{-\int_t^s c(r, X(r ; t, x)) \rmd r} \rmd s\\
    &+g(X(T ; t, x)) e^{-\int_t^T c(r, X(r ; t, x)) \rmd r},
  \end{aligned}
\end{align}
where $X(\cdot) \equiv X(\cdot; t, x)$ is a strong solution of the following stochastic differential equation:
\begin{align}
  \begin{aligned}\label{eq:sde}
    \rmd X(s)&=b(s, X(s)) \rmd s+a(s, X(s)) \rmd W(s), \quad s \in[t, T],\\
    X(t)&=x,      
  \end{aligned}
\end{align}
where $(t, x) \in[0, T) \times \mathbb{R}^d$ and $W(\cdot)$ is a standard $m$-dimensional Wiener process starting at $W(t)=0$ at time $t$.

As stated below, the Feynman--Kac formula guarantees that $\widehat u$ appropriately approximates the solution $v$ of the PDE.
For any $(t,x)\in[0,T]\times\bbR^d$, let 
  $\widehat{\xi}(t,x)\coloneqq\widehat{u}(t,x)-v(t,x)$:
\begin{align}
  \begin{aligned}\label{eq:Feynman--Kac}
  \widehat u(t,x) = v(t, x) + \widehat \xi(t,x),\quad (t, x) \in[0, T] \times \mathbb{R}^d,
  \end{aligned}
\end{align}
where $\widehat{u}$ is defined as \cref{eq:FK_one_sample} and $v$ satisfies \cref{eq:BPDE,eq:BPDE-terminal}.

\begin{proposition}[Feynman--Kac formula for Kolmogorov PDE \cite{Yong1999Stochastic}]
  Suppose that \cref{Asmp:pde} holds.
  Then, $\widehat \xi(t,x)$ satisfies
  \begin{align}\label{eq:FK_noise_mean}
    \bbE[\widehat \xi(t,x)|t,x,v] = 0.
  \end{align}
\end{proposition}

In this study, numerical computation of \cref{eq:FK_one_sample} is called \emph{Feynman--Kac (FK) sampling}.
Note that in the implementation of FK sampling, the computation of \cref{eq:sde} requires a numerical calculation of the SDE using a method such as the Euler--Maruyama method, and the computation of \cref{eq:FK_one_sample} requires a numerical integration using a technique such as the trapezoidal rule.

We describe how to use the data resulting from FK sampling as a dataset $\overline U$ and noise variance matrix $r_{XX}$ that satisfy \cref{eq:observed_data_set}.
For the data set $\overline U$ in \cref{eq:mean_data_set}, we use ${\widecheck U}=[{\widecheck u} (x_1),\ldots,{\widecheck u} (x_N)]^\top$, where ${\widecheck u} (x_i)$ is the sample mean of $M$ FK samples $\{\widehat u_j(0,x_i)\}_{j=1}^M$:
\begin{align}\label{eq:our_estimates}
  \widecheck u{(x_i)}\coloneqq\sum_{j=1}^M \frac{{\widehat u}_j (0,x_i)}{M}, \quad i\in\{1, \ldots, N\},
\end{align}
where $\widehat u_j$ denotes the $j$th FK sample.
It is also appropriate to use $\widecheck r_{XX}$ for the noise variance matrix $r_{XX}$ in \cref{eq:observed_data_set}, which satisfies $[\widecheck r_{XX}]_{ij}=\widecheck r(x_i)\delta_{ij}$, where
\begin{align}\label{eq:our_noise_variance}
  \widecheck r(x_i) \coloneqq \sum_{j=1}^M \frac{(\widehat u_j(0,x_i) - \widecheck u(x_i))^2}{M-1}, \quad i\in\{1, \ldots, N\},
\end{align}
is the unbiased sample variance of $M$ FK samples.
In summary, our dataset satisfies the following:
\begin{align}\label{eq:FK_mean_regression_model}
  \widecheck u(x_i) = v(0,x_i) + \widecheck \xi_i, \quad i\in\{1, \ldots, N\},
\end{align}
where $\widecheck \xi_i\coloneqq \sum_{j=1}^M \widehat \xi_j(0,x_i)/M$ and $\widehat \xi_j(0,x_i) \coloneqq \widehat u_j(0,x_i) - v(0,x_i)$.

\begin{remark}
  Because the noise term $\widecheck \xi$ in \cref{eq:FK_mean_regression_model} is not normally distributed, \cref{eq:mean_regression_model} and \cref{eq:FK_mean_regression_model} do not correspond directly.
  However, if the variance of $\widehat \xi$ is bounded and the sample size $M$ is large enough, the central limit theorem guarantees that the random variable $\widecheck \xi$ approximately follows a normal distribution with mean $0$ and variance $\widecheck r/M$~\cite{kallenberg1997foundations}.
  Because of this fact, we use \cref{eq:FK_mean_regression_model} as an approximate model for \cref{eq:mean_regression_model}.
  In addition, several results assess the non-normality of variables for finite sample sizes $M$~\cite{kallenberg1997foundations,Li2017General}.
  For example, the error of the distributions for finite $M$ is estimated by the Berry--Esseen theorem~\cite{kallenberg1997foundations}, which states that the rate of convergence is $O(M^{-1/2})$.  
\end{remark}

\section{Posterior Variance Analysis}\label{sec:theory}

In this section, we carry out an a priori evaluation of the posterior variance in \cref{eq:posterior_variance}.
We focus on the integral value of the MSE on the input space $\calX$, called integrated MSE (IMSE):
\begin{align}\label{eq:IMSE}
  \mathrm{IMSE}\coloneqq\int_{\calX} \widetilde\sigma^2(x) \rmd \nu(x),
\end{align}
where $\nu$ is a given probability measure of the observation points whose support is $\calX$.

Here, we aim to characterize the speed of the decay of IMSE with respect to an increase in the number of observation points $N$ and the sample size $M$.
To achieve this, we first characterize the decay speed using the eigenvalues of the kernel function and the magnitude of the noise term.
Specifically, we derive a probabilistic lower bound for the IMSE and show that the IMSE is expected to decrease as $N$ and $M$ increase (see \cref{thm:convergence-MSE-hetero-min}).
Among the variables used in this lower bound, the eigenvalues can be considered known, but the magnitude of the noise term is unknown.
We therefore probabilistically assess the error between the true value of the noise term and its estimate, and show that the evaluation of the IMSE is justified under sufficiently large sample size $M$ (see \cref{thm:err-bound}).
After deriving these theorems, we give a more conservative estimate of the IMSE, which does not include information on the eigenvalues of the kernel and the magnitude of the noise term (see \cref{col:conservative-bound}).

\begin{assumption}\label{assmp:convergence-MSE-hetero-min}
  \begin{enumerate}
    \item The kernel function $k$ is bounded and positive-definite on $\calX$.
    \item The observation points $X=[x_1,\ldots,x_N]^\top$ are sampled from the known probability measure $\nu$.
    \item The data set $\overline U$ in \cref{eq:mean_data_set} satisfies \cref{eq:observed_data_set}.
    \item For the noise variance function $r$ in \cref{eq:noise-var-hetero}, there exists a positive minimum value $r_\mathrm{min}=\min_{x\in\calX} r(x)>0$. 
  \end{enumerate}
\end{assumption}

\begin{theorem}\label{thm:convergence-MSE-hetero-min}
  Suppose that \cref{assmp:convergence-MSE-hetero-min} holds.
  Then, for any $\epsilon>0$, there exists an integer $\tilN(\epsilon)>0$ such that for any $N>\tilN$ the following holds:
  \begin{align}\label{eq:IMSE-bound-hetero-min}
    \mathrm{Pr} \left( \mathrm{IMSE}  > (1-\epsilon) L_\mathrm{IMSE}\right)= 1,
  \end{align}
  where 
  \begin{align}\label{eq:IMSE-bound-hetero-min2}
    L_\mathrm{IMSE}\coloneqq r_\mathrm{min} \sum_{p \in I} \frac{\lambda_p}{r_\mathrm{min} + NM \lambda_p}.
  \end{align}
  Here, $\{\lambda_p\}_{p\in I}$ denotes the eigenvalues of the kernel function $k$ for the probability measure $\nu$.
  The set $I\subseteq \bbN$ denotes the index set of the eigenvalues of the kernel. 
\end{theorem}

\begin{remark}
  In ordinary error analysis in the numerical method of PDEs, upper bounds on the decay of the error are derived. 
  The lower bound given here helps assess the possible best performance of the method. 
  For example, one can estimate the minimum number of observation points needed for a pre-determined performance requirement.
\end{remark}

\begin{remark}
  The evaluation of lower bounds on the posterior variance has long been studied as an asymptotic analysis of Gaussian processes~\cite{Sollich2002Learning,williams2006gaussian}.
  To the best of our knowledge, they only address the case where the noise model in the GPR is uniform.
  We extend the results to the heteroscedastic case.
\end{remark}

\begin{remark}
  In the right-hand side of \eqref{eq:IMSE-bound-hetero-min}, the actual decay speed depends on the decay rate of the eigenvalues of the kernel. 
  For some smooth kernels, the decay rate of the eigenvalues has been evaluated under known measures of the observation points~\cite{Santin2016Approximation}. 
  It is also possible to evaluate the eigenvalues numerically by performing a sampling routine (see chapter 4 in \cite{williams2006gaussian}).
\end{remark}

To evaluate the right-hand side of \cref{thm:convergence-MSE-hetero-min}, it is necessary to find the minimum 
$r_{\min}$
of the noise variance function $r$ of the GPR.
However, as mentioned in \cref{sec:FK}, $r$ can only be estimated pointwise using sample variance as in \cref{eq:our_noise_variance}.
Therefore, it is natural to use 
$\overline r_{\min} \coloneqq \overline {r}(x_{i^*})\ (i^* \coloneqq \argmin_{i\in\{1,\ldots,N\}}\overline {r}(x_i))$ as an estimate of $r_{\min}$, where $\overline r(x_i) = \sum_{j=1}^M (u_j(x_i)-\overline u(x_i))^2/(M-1)$ is the unbiased variance of \cref{eq:regression_model}.
The following theorem probabilistically evaluates the error between $r_{\min}$ and $\overline r_{\min}$.

\begin{theorem}\label{thm:err-bound}
  Suppose that \cref{assmp:convergence-MSE-hetero-min}.2--\ref{assmp:convergence-MSE-hetero-min}.4 holds. 
  Then, for any $N$, $M$, $\delta>0$, and $\epsilon>\delta$,
  \begin{align}\label{eq:err-bound}
    \begin{aligned}  
      &\mathrm{Pr}(|\overline r_{\min} - r_{\min}|\ge \epsilon, \delta>|r(x_{i^*}) - r_{\min}|)
      \\
      &\le 1-\prod_{i=1}^N \max \left\{1-\frac{2 r^2\left(x_i\right)}{(\epsilon-\delta)^2(M-1)}, 0\right\}.
    \end{aligned}
  \end{align}
\end{theorem}

\begin{remark}
  From \cref{eq:err-bound}, we can justify replacing $r_{\min}$ in \cref{eq:IMSE-bound-hetero-min2} with the minimum sample variance $\overline r_{\min}$ under a sufficiently large $M$ relative to $N$.
\end{remark}

Apart from the results above, it is also possible to evaluate the lower bound of the posterior variance in a form that does not include eigenvalues or noise variance.
The following lower bound is more conservative than the one given in \cref{thm:convergence-MSE-hetero-min} but is more useful for an intuitive understanding.

\begin{corollary}\label{col:conservative-bound}
  Suppose that \cref{assmp:convergence-MSE-hetero-min} holds.
  Then, there exists a constant $C>0$ such that for any $N\in\bbN$ and $M\in\bbN$, 
  \begin{align}\label{eq:IMSE-lowerbound-hetero-conservative}
    L_\mathrm{IMSE} \ge \frac{C}{NM}.
  \end{align}
\end{corollary} 

\begin{remark}
  From \cref{eq:IMSE-lowerbound-hetero-conservative}, we see that the fastest convergence rate of the posterior variance is $O((NM)^{-1})$ for an FK sample size $M$ and the number of GPR observation points $N$.
\end{remark}

\section{Numerical Calculation}\label{sec:numerics}

In this section, we apply the proposed method to several PDEs and evaluate their performance.
The target PDEs are the heat equation, the advection-diffusion equation, and the Hamilton--Jacobi--Bellman (HJB) equation.
The goal is to learn the functional form of the solution on $x\in[0,1]^d$ for $d=10$ spatial dimensions.
Specifically, we learn only the first-dimensional component of the PDE solution; we use $N$ equally spaced grid points in the range $[0,1]$ for the $1$-dimensional component and fix the $2$--$9$-dimensional components to $0.5$.
In other words, we set $\calX=[0, 1]\times{0.5}\times\cdots\times{0.5}$ and we use a one-dimensional Lebesgue measure as the probability measure $\nu$.
This is because using $\widetilde N = N^d$ as the number of observation points would drain the machine's memory for allocating variables in GPR.
The approximation methods proposed in \cite{Gardner2018Producta} would solve this problem, but this implementation is beyond the scope of this paper and is not done here.

We compare the performance of the following three estimation methods:
\begin{itemize}
  \item Heteroscedastic GPR (HSGPR): The proposed method, where the noise variance function in GPR is estimated using the sample variance of the FK sample, using \cref{eq:our_noise_variance}.
  \item Standard GPR: A method that treats $r\equiv\sigma^2\ (\sigma\in\bbR)$ as a hyperparameter in GPR, where the hyperparameters are determined to maximize the marginal likelihood.
  \item Linear interpolation: A method to linearly interpolate the function value $\{\overline u(x_i)\}_{i=1}^N$ obtained by FK sampling.
\end{itemize}
We evaluate each method using the following criteria:
\begin{enumerate}
  \item Squared $L^2$-norm error of the estimates for the test data:
  \begin{align}\label{eq:l2-error}
    \mathrm{Error} = {\int_\calX \qty[\tilu(x) - v(x)]^2 \rmd x},
  \end{align}
  where $\tilu(x)$ is the estimate of the solution of the PDE and $v(x)$ is the test data to approximate the true solution $v(0,x)$ of the PDE \cref{eq:BPDE,eq:BPDE-terminal}.
  For numerical integration in \cref{eq:l2-error}, the trapezoidal law is used.

  \item Integrated mean squared error (IMSE) defined in \cref{eq:IMSE}, along with its bound in \cref{eq:IMSE-bound-hetero-min2}.
  Here, to evaluate the bound in \cref{eq:IMSE-bound-hetero-min2}, the eigenvalues are approximated using a sampling routine, as described in chapter 4 in \cite{williams2006gaussian}.
  The term $r_{\min}$ is approximated with $\widecheck r_{\min}=\min_{i} \widecheck r(x_i)$, where $\{\widecheck r(x_i)\}_{i=1}^N$ are the sample variance of the FK sample \cref{eq:FK_one_sample}.
\end{enumerate}
The first criterion evaluates the accuracy of the estimate, while the second criterion evaluates the uncertainty.
Here, the results of FK sampling with sufficiently large $M,N$, i.e. $M=40000$ and $N=40$, are used as test data $v$.
Note that only GPR-based methods can calculate the IMSE.

In GPR, we use the Mat\'{e}rn kernel corresponding to a Sobolev space of order $2$, considering that we seek to solve second-order PDEs.
The hyperparameters in the kernel are determined to maximize the marginal likelihood of the data.
This maximization is performed by a gradient ascent method on the log-likelihood of the hyperparameters.
Since the log-likelihood is a non-convex function, the optimization steps converge to a local solution depending on the choice of initial values.
To search for a better local solution, $27$ initial values of hyperparameters are randomly determined, and the value with the largest likelihood is adopted.
Functions built into {\ttfamily gpytorch} are used for optimizing hyperparameters, calculating posterior distributions and kernel functions~\cite{Gardner2018GPyTorch}.

In FK sampling, the Euler--Maruyama method is used to numerically calculate the SDE of \cref{eq:sde} at each observation point, and the trapezoidal law is used for numerical integration in \cref{eq:Feynman--Kac}.
The results of the FK sampling depend on the sample path of the Brownian motion.
To suppress the effect of this randomness, we ran $50$ simulations with different random number seeds in each experiment and compared the results using the average of these simulations.

Before discussing the results for each PDE, we give a summary of the results that are common to each experiment.
First, regarding the squared $L^2$-norm error, the error tends to become smaller as $M, N$ are increased for all three methods.
The errors for HSGPR and GPR are comparable, and these two methods achieve smaller errors than linear interpolation.
With respect to IMSE, both HSGPR and GP achieve smaller IMSE for larger $N$ and $M$.
Comparing the two, HSGPR achieves a smaller IMSE overall, suggesting that the proposed method achieves estimation with lower uncertainty.

\subsection{Heat Equation}
  We considered the following initial value problem for the heat equation:
  \begin{align}
    &\begin{aligned}\label{eq:heat}
      \partial_t w(t,x)&=\frac{1}{2} \operatorname{tr}\left\{a a^{\top} \partial_{xx} w(t, x)\right\},
    \end{aligned} \quad (t, x) \in[0, T) \times \mathbb{R}^d, \\
    &w(0,x)=\exp \left(-\lambda(x-0.5)^{\top} (x-0.5)\right), \qquad\qquad x \in\mathbb{R}^d,\label{eq:heat-initial}
  \end{align}
  where we have defined $a \equiv 0.4 I_d$, $\lambda= 5.0$, and $T=1.0$.
  To solve the initial value problem, the proposed method is applied after using the variable transformation described in \cref{rem:initial-terminal}.

  \begin{figure}[ht]
    \centering
    \begin{subfigure}{0.45\linewidth}
    \centering
      \includegraphics[width=1\linewidth]{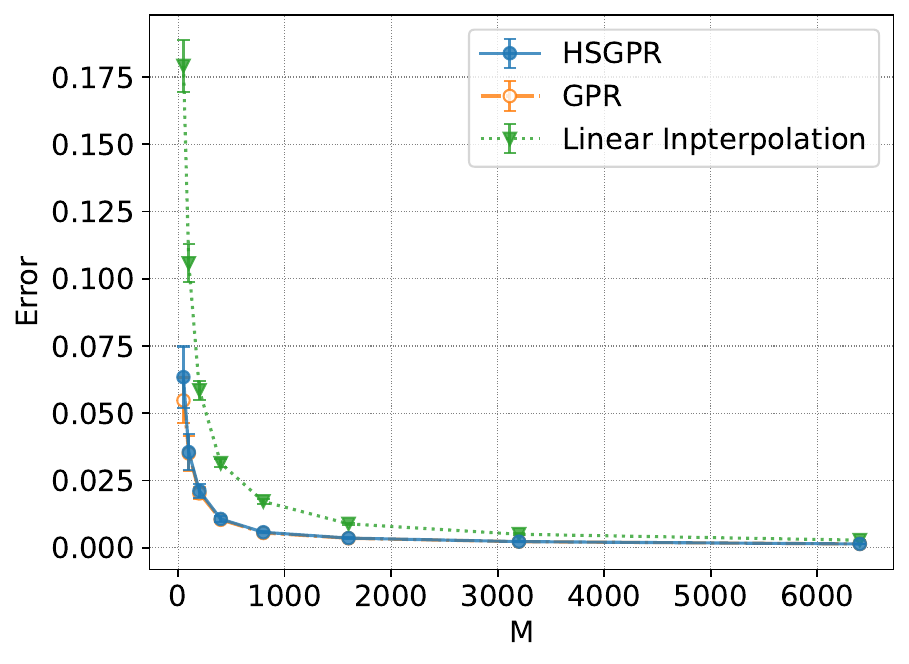}
    \caption{Squared $L^2$-norm error for the Feynman--Kac sample size $M$.}
    \end{subfigure}
    \begin{subfigure}{0.45\linewidth}
    \centering
      \includegraphics[width=1\linewidth]{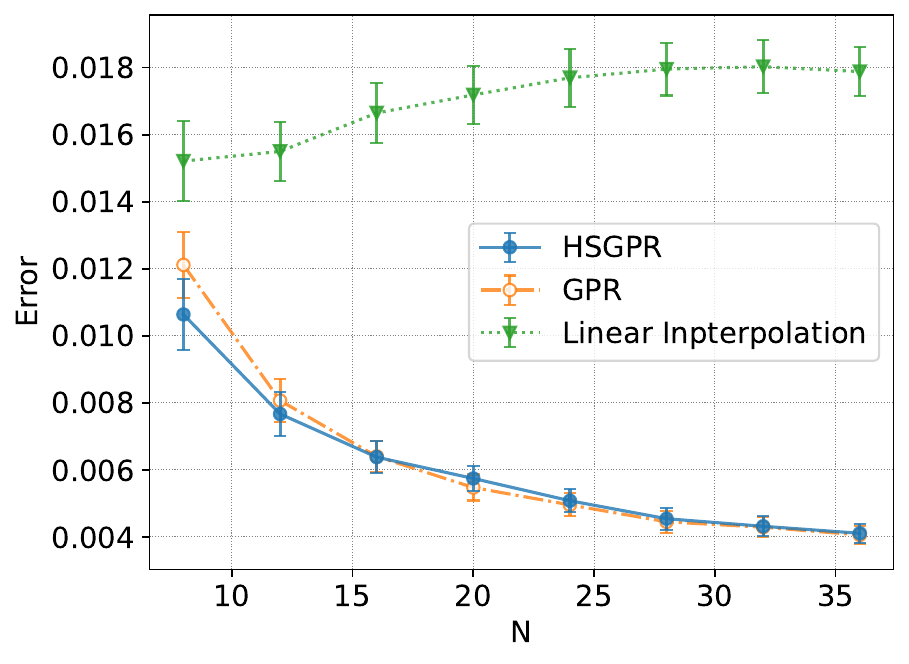}
    \caption{Squared $L^2$-norm error for the number of observation points $N$.}
    \end{subfigure}
    \begin{subfigure}{0.45\linewidth}
      \centering
        \includegraphics[width=1\linewidth]{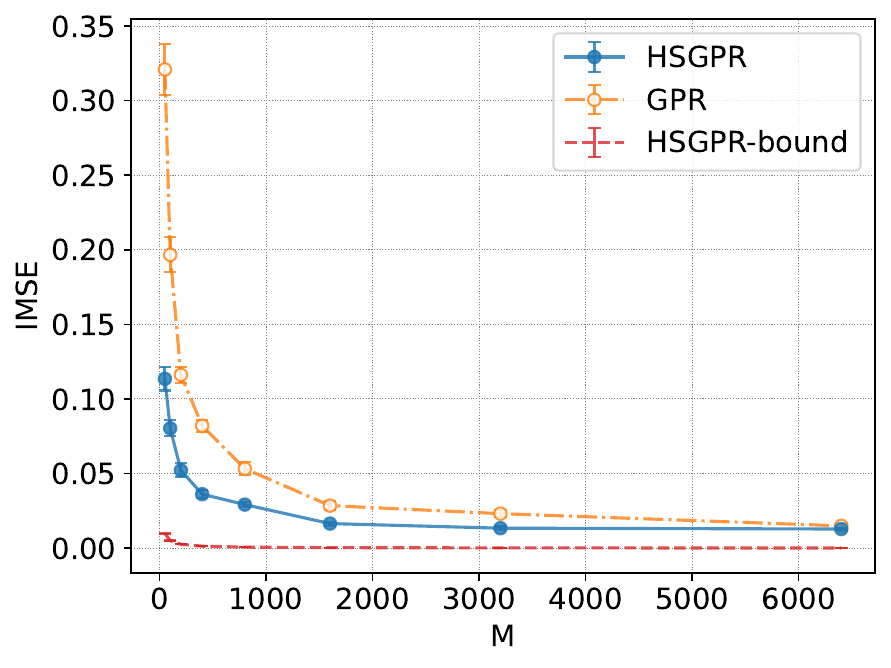}
      \caption{IMSE and its bound for the Feynman--Kac sample size $M$.}
      \end{subfigure}
      \begin{subfigure}{0.45\linewidth}
      \centering
        \includegraphics[width=1\linewidth]{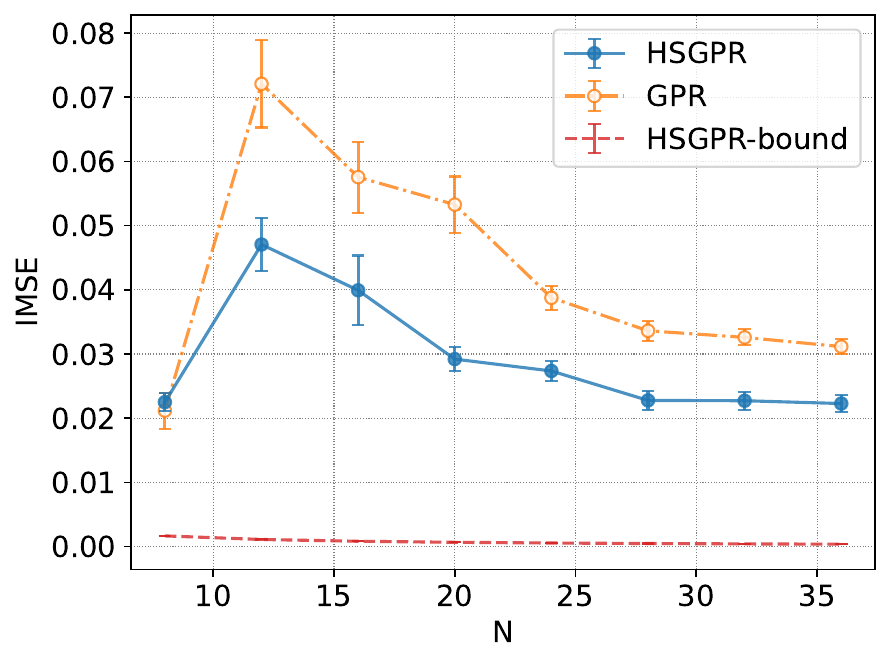}
      \caption{IMSE and its bound for the number of observation points $N$.}
      \end{subfigure}
    \caption{Squared $L^2$-norm error and IMSE when each method is applied to the 10-dimensional heat equation. Blue (solid): HSGPR / Orange (dash-dot): GPR / Green (dotted): linear interpolation / Red (dashed): lower bound of IMSE in HSGPR.
    The error bars represent the standard error of the calculation results for 50 different seeds used in the FK sampling.
    When varying the sample size $M$, the number of observation points $N$ is set as $N=20$, and when varying the number of observation points $N$, the sample size $M$ is set as $M=800$.
    }
    \label{fig:heat_1d_err_likelihood}
  \end{figure}

  We examined the behavior of the squared $L^2$-norm error and IMSE with respect to sample size $M$ and number of observation points $N$. 
  Here, we set $N=20$ when varying $M$ and $M=800$ when varying $N$. 
  The respective results are shown in \cref{fig:heat_1d_err_likelihood}.
  We see that increasing $M$ or $N$ in each method tends to decrease the $L^2$-norm error and IMSE.

  Regarding the $L^2$-norm error, all methods decrease the error as the sample size $M$ is increased.
  This shows that the proposed method approximates the solution of the PDE more accurately as the parameters are increased.
  The errors for HSGPR and GPR are comparable, while the errors for linear interpolation are worse.
  As the number of observation points $N$ is increased, linear interpolation tends to have a larger error, while HSGPR and GPR have smaller errors.
  In the region where $N$ is small, HSGPR reduces the error more than GPR.
  
  For both HSGPR and GPR, the IMSE decreases monotonically with increasing $M$.
  For all $M$, HSGPR achieves a smaller IMSE than GPR. 
  When $N$ is increased, the IMSE decreases in both HSGPR and GPR, and the value in HSGPR is smaller than the value in GPR.
  The approximated IMSE bounds are found to decay on the same order of magnitude as the measured values.
  
  Summarizing the results for both criteria, HSGPR is better than linear interpolation but similar to GPR in terms of error, and better than GPR in terms of uncertainty.

\subsection{Advection-Diffusion Equation}
  
  Next, we consider the following initial value problem for the advection-diffusion equation:
  \begin{align}
    &\begin{aligned}\label{eq:fp}
      \partial_t w(t,x)&=\frac{1}{2} \operatorname{tr}\left\{a a^{\top} \partial_{xx} w(t, x)\right\} + b^\top w(t, x),
    \end{aligned} \quad (t, x) \in[0, T) \times \mathbb{R}^d, \\
    &w(0,x)=\exp \left(-\lambda(x-0.5)^{\top} (x-0.5)\right), \qquad\qquad\quad x \in\mathbb{R}^d,\label{eq:fp-initial}
  \end{align}
  where we have defined $a \equiv 0.4 I_d$ and $b=0.01\mathbbm{1}_d$, $\lambda=5.0$, and $T=1.0$.
  As in the previous subsection, the proposed method is applied after using the variable transformation described in \cref{rem:initial-terminal}.

  \begin{figure}[ht]
    \centering
    \begin{subfigure}{0.45\linewidth}
    \centering
      \includegraphics[width=1\linewidth]{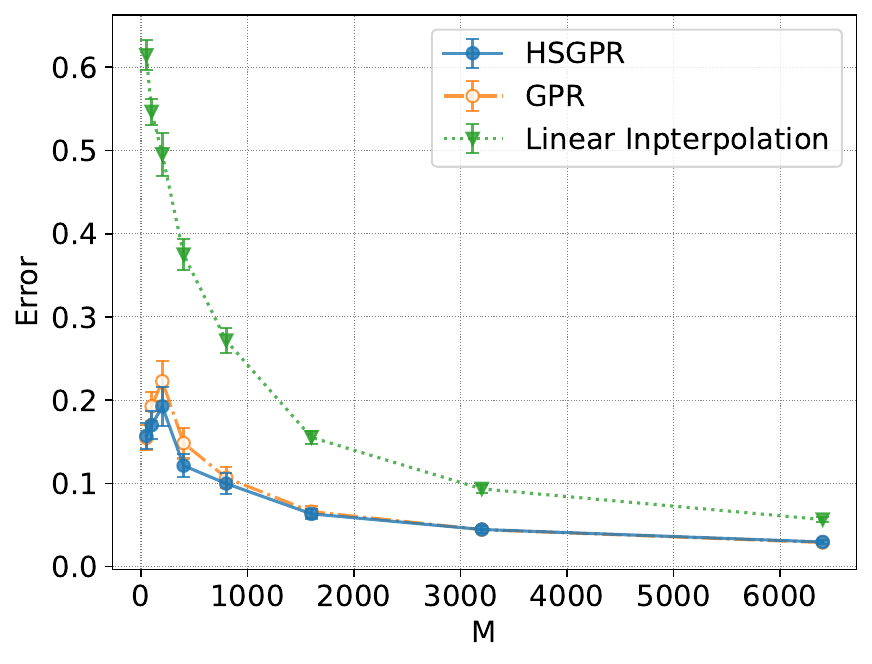}
    \caption{Squared $L^2$-norm error for the Feynman--Kac sample size $M$.}
    \end{subfigure}
    \begin{subfigure}{0.45\linewidth}
    \centering
      \includegraphics[width=1\linewidth]{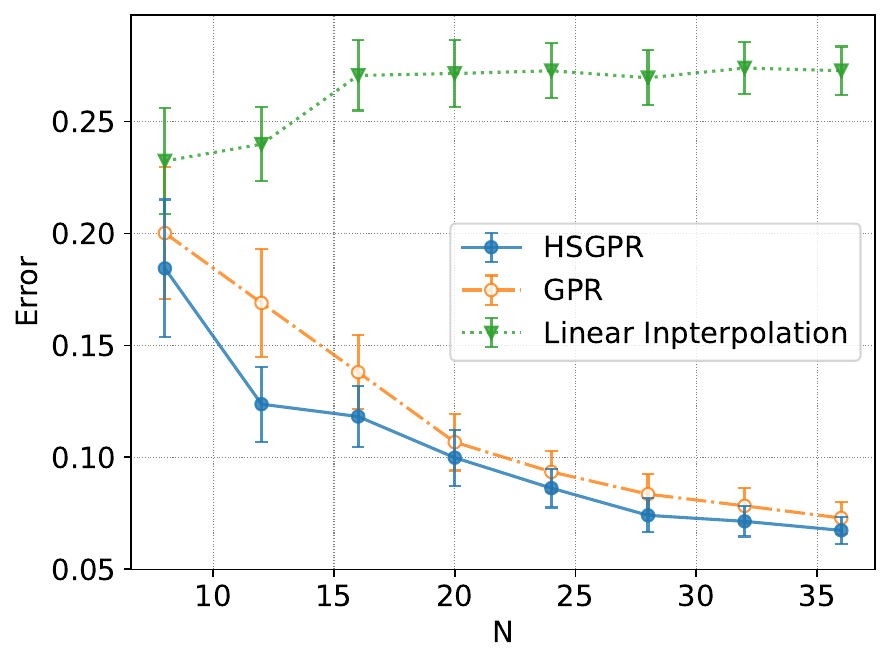}
    \caption{Squared $L^2$-norm error for the number of observation points $N$.}
    \end{subfigure}
    \begin{subfigure}{0.45\linewidth}
      \centering
        \includegraphics[width=1\linewidth]{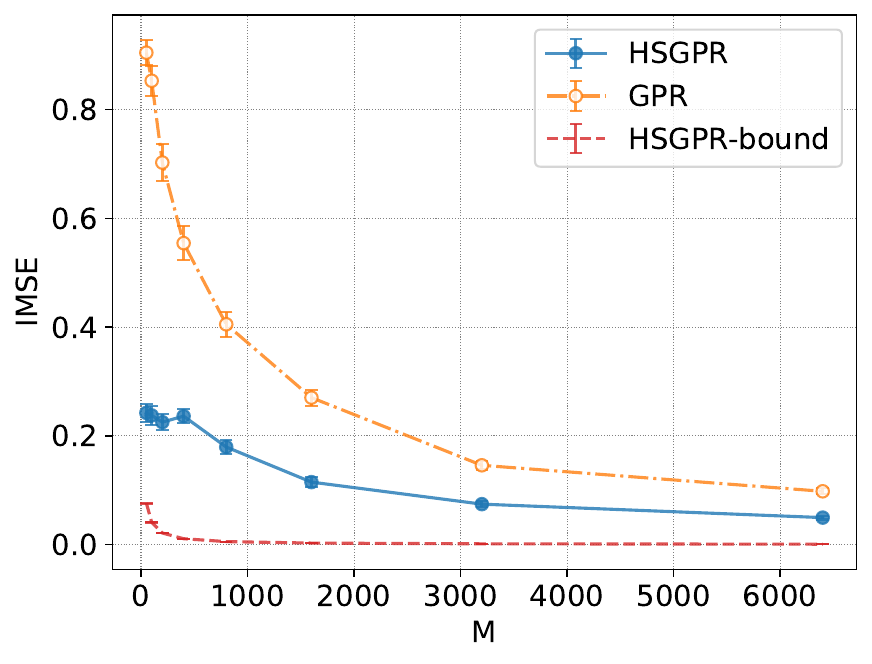}
      \caption{IMSE and its bound for the Feynman--Kac sample size $M$.}
      \end{subfigure}
      \begin{subfigure}{0.45\linewidth}
      \centering
        \includegraphics[width=1\linewidth]{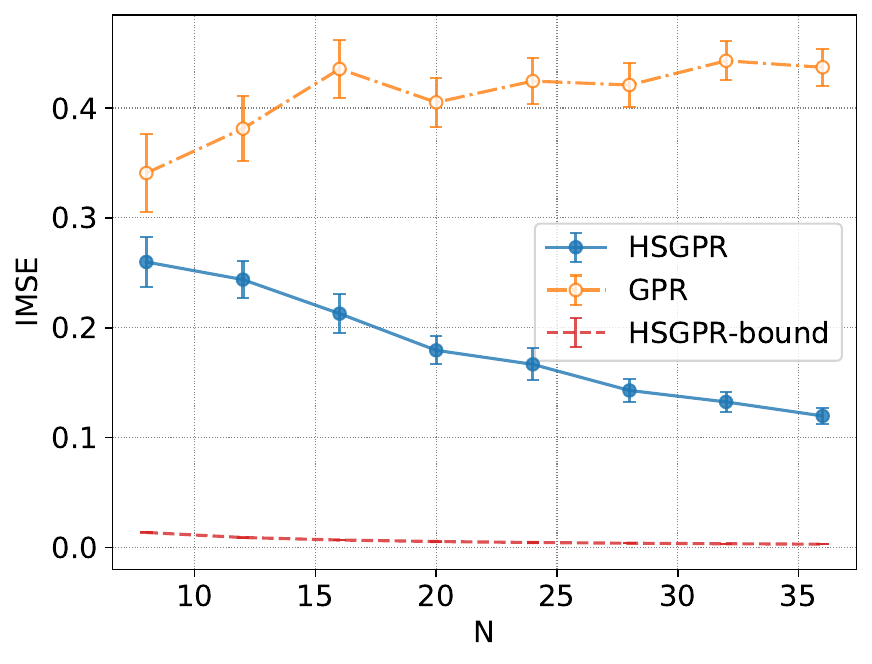}
      \caption{IMSE and its bound for the number of observation points $N$.}
      \end{subfigure}
    \caption{Squared $L^2$-norm error and IMSE when each method is applied to the 10-dimensional advection-diffusion equation. Blue (solid): HSGPR / Orange (dash-dot): GPR / Green (dotted): linear interpolation / Red (dashed): lower bound of IMSE in HSGPR.
    The error bars represent the standard error of the calculation results for 50 different seeds used in the FK sampling.
    When varying the sample size $M$, we set the number of observation points $N$ as $N=20$, and when varying the number of observation points $N$, we set the sample size $M$ as $M=800$.}
    \label{fig:fp_1d_err_likelihood}
  \end{figure}

  The squared $L^2$-norm error and IMSE when varying $M$ and $N$ are plotted in \cref{fig:fp_1d_err_likelihood}.
  Similar to the example above, increasing $M$ or $N$ in each method tends to decrease the $L^2$-norm error and IMSE.

  All methods decrease the $L^2$-norm error as the sample size $M$ is increased.
  The errors for HSGPR are slightly smaller than for GPR, while the errors for linear interpolation are larger.
  As the number of observation points $N$ is increased, linear interpolation tends to have a larger error, while HSGPR and GPR have smaller errors.
  HSGPR reduces the error more than GPR, particularly in the region where $N$ is small. 
  
  For both HSGPR and GPR, the IMSE decreases monotonically when increasing $M$, and for all $M$, HSGPR achieves a smaller IMSE than GPR. 
  When $N$ is increased, the IMSE of HSGPR decreases, while that of GPR increases.
  The approximately calculated IMSE bounds decay on a similar order of magnitude to the measured values.

\subsection{Hamilton--Jacobi--Bellman Equation}

  Finally, we considered the following terminal value problem for the Hamilton--Jacobi--Bellman (HJB) equation:
  \begin{align}
      -\partial_t v(t, x) & = \ell(x) -\frac{1}{2} \partial_x v(t, x)^{\top} B R^{-1} B^{\top} \partial_x v(t, x)+\operatorname{tr}\left\{aa^\top \partial_{x x} v(t, x)\right\},\nonumber\\
      &\qquad\qquad\qquad\qquad\qquad\qquad\qquad\qquad\quad (t, x) \in[0, T) \times \mathbb{R}^d, 
    \label{eq:hjb-0}\\
    v(T, x) & =g(x),\qquad\qquad\qquad\qquad\qquad\qquad\qquad x \in\mathbb{R}^d, \label{eq:hjb-0-initial}
  \end{align}
  where we have defined $\ell(x) = (x-0.5\mathbbm{1}_d)^\top (x-0.5\mathbbm{1}_d)$, $g(x)\equiv 0$, $B=I_d$, $R=I_d$, $a \equiv 0.4 I_d$, and $T=1.0$.
  This PDE corresponds to the following optimal control problem:
  \begin{align}
    &\text{minimize } J(\alpha) = \mathbb E\left\{\int_{0}^{T} \ell(X(s)) + \alpha(s)^\top R \alpha(s) \mathrm{d} s + g (X(T))\right\},\\
    &\text{subject to } \rmd X (s)=B\alpha(s)\rmd s + a \rmd W(s), \ \ s \in[0, T],\\
    &\qquad\qquad\quad X(0) = x.
  \end{align}
  Since the HJB equation contains nonlinear terms, the proposed method is not directly applicable.
  Thus, we use the Cole--Hopf transformation to linearize the equation.
  First, we find the real number $\lambda$ which satisfies $\lambda BR^{-1}B^\top = a a^\top$.
  Using this, we transform the PDE with $\tilv=\exp\qty(-\frac{v}{\lambda})$:
  \begin{align}
    &\begin{aligned}\label{eq:hjb}
      \partial_t \tilv(t,x)&=-\frac{1}{2} \operatorname{tr}\left\{a a^{\top} \partial_{xx} \tilv(t, x)\right\} + \frac{1}{\lambda} \ell(x) \tilv(t,x),
    \end{aligned}  (t, x) \in[0, T) \times \mathbb{R}^d, \\
    &\tilv(T,x)=\exp\qty(-\frac{g(x)}{\lambda}), \qquad\qquad\qquad\qquad\qquad\qquad x \in\mathbb{R}^d,\label{eq:hjb-initial}
  \end{align}
  We conducted FK sampling on this linearized HJB equation and performed the regression for the Cole-Hopf inverse transformed data.

  \begin{figure}[ht]
    \centering
    \begin{subfigure}{0.45\linewidth}
    \centering
      \includegraphics[width=1\linewidth]{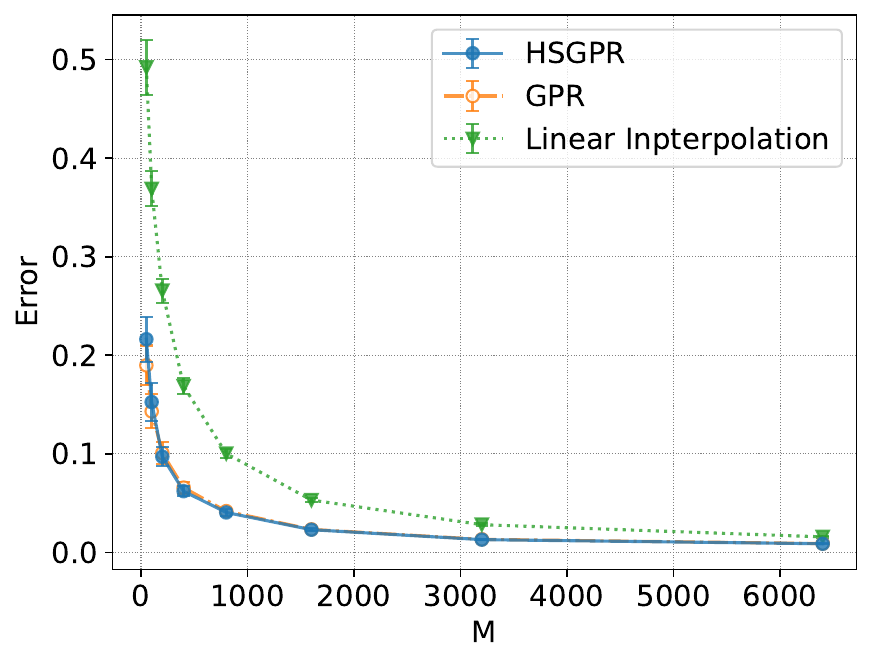}
    \caption{Squared $L^2$-norm error for the Feynman--Kac sample size $M$.}
    \end{subfigure}
    \begin{subfigure}{0.45\linewidth}
    \centering
      \includegraphics[width=1\linewidth]{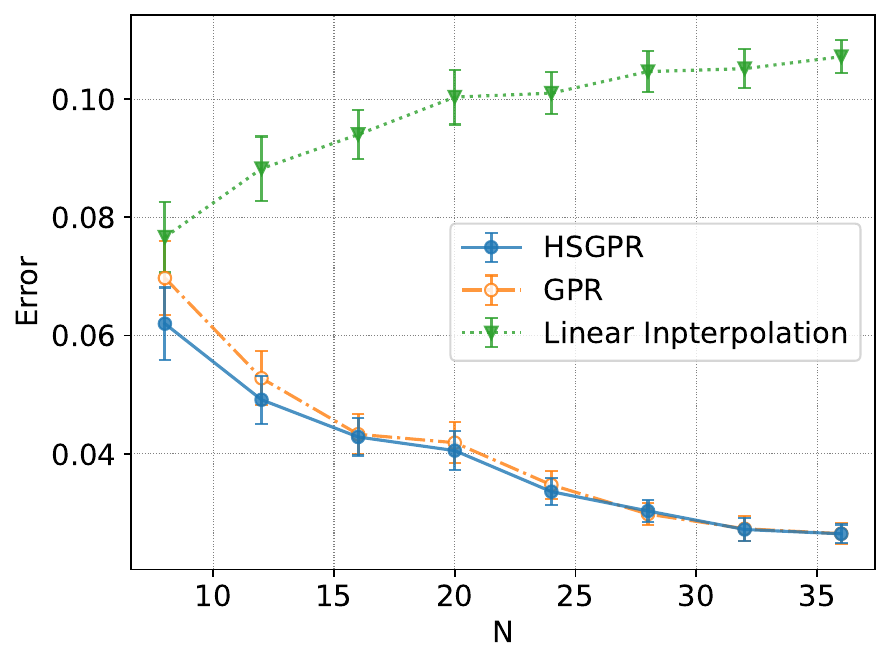}
    \caption{Squared $L^2$-norm error for the number of observation points $N$.}
    \end{subfigure}
    \begin{subfigure}{0.45\linewidth}
      \centering
        \includegraphics[width=1\linewidth]{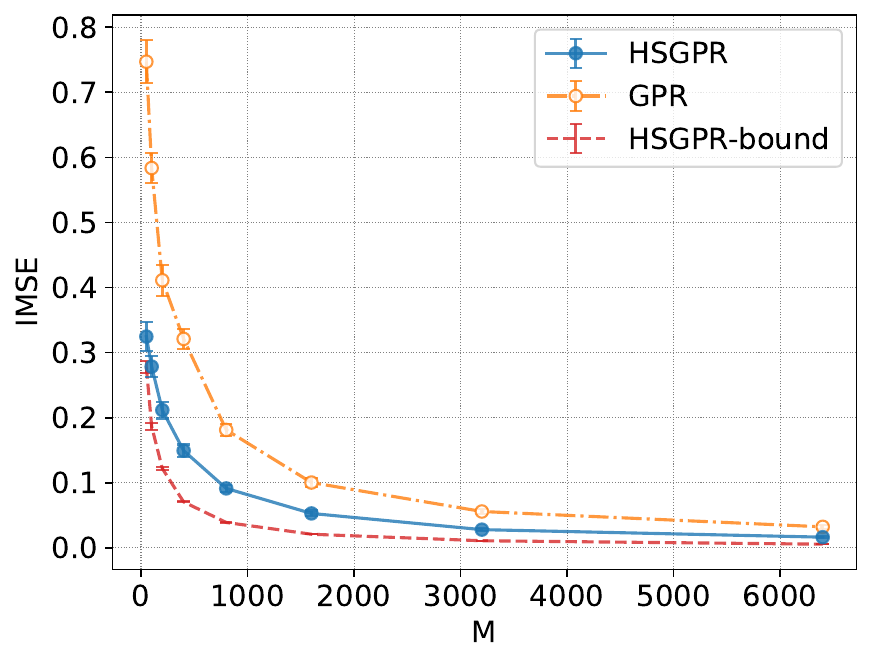}
      \caption{IMSE and its bound for the Feynman--Kac sample size $M$.}
      \end{subfigure}
      \begin{subfigure}{0.45\linewidth}
      \centering
        \includegraphics[width=1\linewidth]{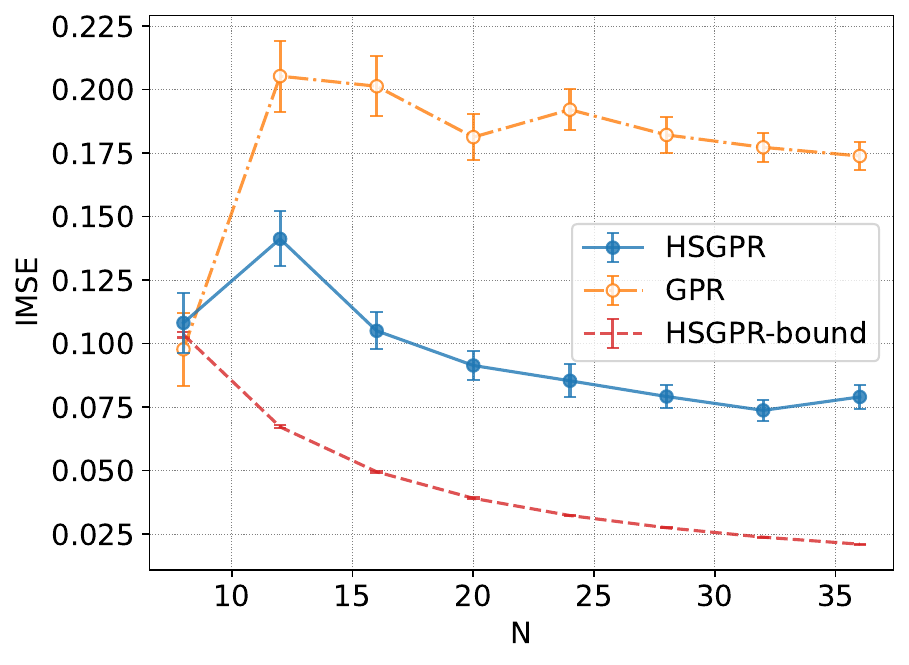}
      \caption{IMSE and its bound for the number of observation points $N$.}
      \end{subfigure}
    \caption{Squared $L^2$-norm error and IMSE when each method is applied to the 10-dimensional HJB equation. Blue (solid): HSGPR / Orange (dash-dot): GPR / Green (dotted): linear interpolation / Red (dashed): lower bound of IMSE in HSGPR.
    The error bars represent the standard error of the calculation results for 50 different seeds used in the FK sampling.
    When varying the sample size $M$, we set the number of observation points $N$ as $N=20$, and when varying the number of observation points $N$, we set the sample size $M$ as $M=800$.}
    \label{fig:hjb_1d_err_likelihood}
  \end{figure}

  The squared $L^2$-norm error and IMSE when varying $M$ and $N$ are plotted in \cref{fig:hjb_1d_err_likelihood}.
  Similar to the two examples above, increasing $M$ or $N$ in each method tends to decrease the $L^2$-norm error and IMSE.

  All methods decrease the $L^2$-norm error as sample size $M$ is increased.
  The values for HSGPR and GPR are comparable, while the values for linear interpolation are larger.
  When the number of observation points $N$ is increased, linear interpolation has a larger error, while HSGPR and GPR have smaller errors.
  HSGPR reduces the error more than GPR, particularly in the region where $N$ is small. 
  
  Next, for both HSGPR and GPR, the IMSE decreases monotonically with increasing $M$, and for all $M$, HSGPR achieves a smaller IMSE than GPR. 
  As $N$ is increased, the IMSE of the HSGPR decreases more than that of the GPR.
  The tightness of IMSE bounds is improved compared to the two examples above.
  In the region where $N$ is small, the IMSE of the GPR is close to the estimated lower bound.
  Note that the lower bound given in \cref{thm:convergence-MSE-hetero-min} is a probabilistic inequality that holds when $N$ is sufficiently large and does not necessarily hold for small $N$.

\section{Conclusion}\label{sec:fkgp-conclusion}

We have proposed a numerical method for PDEs by combining the heteroscedastic Gaussian process regression and the Feynman--Kac formula.
The proposed method is able to assess the uncertainty in the numerical results and to seek high-dimensional PDEs.
The quality of the solution can be improved by adjusting the kernel function and incorporating noise information obtained from Monte Carlo samples into the GPR noise model.
The method succeeds in providing a lower bound on the posterior variance, which is useful for estimating the best-case performance in numerical calculations. 
The numerical results confirm that the proposed method provides better approximations in the wide parameter domains compared to naive linear interpolation and standard GPR-based methods.

In future work, we will explore iterative observation schemes that leverage uncertainty information. 
By adaptively sampling the most uncertain points, we expect to efficiently estimate the solution of the PDE. 
We will also investigate the design of uncertainty-aware acquisition functions to guide the search for regions of interest, such as those that maximize or minimize a specific physical quantity.

\section*{Acknowledgment}
We would like to thank Dr.~Yuta Koike of the Graduate School of Mathematical Sciences at the University of Tokyo for helpful comments on the proof of \cref{thm:err-bound}.

\section*{Author contributions}
Conceptualization: D.I.; Methodology: D.I.; Formal analysis and investigation: D.I., Y.I., T.K., and N.S.; Writing - original draft preparation: D.I.; Writing - review and editing: Y.I., T.K, N.S., and H.Y.; Resources: N.S. and H.Y.; Supervision: N.S. and H.Y.

\section*{Declarations}
The authors did not receive support from any organization for the submitted work.
The authors have no competing interests to declare that are relevant to the content of this article.

\section*{Data Availability}
The data that support the findings of this study are available from the corresponding author upon reasonable request.

\section*{Appendix}\label{sec:appendix}


\subsection*{Proof of \cref{thm:convergence-MSE-hetero-min}}
  Our proof is based on the proof of Theorem 1 of \cite{LeGratiet2015Asymptotic}. 
  The main difference is the extension of the noise model in the GPR to the heteroscedastic case.
  Proofs for nondegenerate kernels can be given by way of using the Karhunen--Lo\`eve expansion as in \cite{LeGratiet2015Asymptotic}.
  Thus, we only consider degenerate kernels, that is, the case where the number of nonzero eigenvalues of the kernel is finite.

  Let $\{\lambda_p, \phi_p\}_{p\in I}$ be the eigensystems of the kernel function $k$ for the probability measure $\nu$.
  Let $\overline p$ be the number of nonzero eigenvalues.
  We define $\Lambda\coloneqq\operatorname{diag}\left(\lambda_i\right)_{1 \leq i \leq \overline{p}}$, 
  $\phi(x)\coloneqq\left[\phi_1(x), \ldots, \phi_{\overline{p}}(x)\right]^\top$, and
  $\Phi\coloneqq\left[\phi\left(x_1\right), \ldots, \phi\left(x_{N}\right)\right]^\top$.
  Then, by using \cref{eq:expansion},
  the MSE in \cref{eq:posterior_variance} is rewritten in the following form:
  \begin{align}
    \begin{aligned}
      \widetilde\sigma^2(x) 
      &=\phi(x)^\top \Lambda \phi(x)-\phi(x)^\top \Lambda \Phi^\top\left(\Phi \Lambda \Phi^\top + \frac{r_{XX}}{M}\right)^{-1} \Phi \Lambda \phi(x)\\
      &\ge \phi(x)^\top \Lambda \phi(x)-\phi(x)^\top \Lambda \Phi^\top\left(\Phi \Lambda \Phi^\top + \frac{r_{\min}}{M} I_{N}\right)^{-1} \Phi \Lambda \phi(x),
    \end{aligned}
  \end{align}
  where $I_N\in\bbR^{N\times N}$ denotes the identity matrix of size $N$.
  By using Woodbury-Sherman-Morrison formula~\cite{Petersen2008Matrix}, this is rearranged as
  \begin{align}\label{eq:woodbury}
    \begin{aligned}
      \widetilde\sigma^2(x) 
      &\ge 
      \phi(x)^\top\left( \frac{M}{r_{\min}}\Phi^\top \Phi + \Lambda^{-1}\right)^{-1} \phi(x).
    \end{aligned}
  \end{align}
  On the right-hand side, each component of the matrix $\Phi^\top \Phi$ is written as
  \begin{align}
    (\Phi^\top \Phi)_{p q} = \sum_{i\in\{1,\ldots,N\}} \phi_p(x_i)\phi_q(x_i).
  \end{align}
  By the strong law of large numbers, the right-hand side divided by $N$ converges almost surely to the Kronecker delta (see \cref{eq:eigen-orthogonal}).
  That is, for any $\epsilon> 0$, there exists $\tilN(\epsilon)>0$ and for any $N\ge \tilN(\epsilon)$,
  \begin{align}\label{eq:low-of-large-numbers}
    \mathrm{Pr}\left( \left|\frac{1}{N}\sum_{i\in\{1,\ldots,N\}} \phi_p(x_i)\phi_q(x_i)- \delta_{pq}\right| < \frac{\epsilon}{\overline p} \right) = 1.
  \end{align}
  From this and the Gershgorin circle theorem, we can obtain
  \begin{align}
    \mathrm{Pr}\left( \frac{M}{r_{\min}}\Phi^\top \Phi  + \Lambda^{-1} \prec \frac{NM(1+\epsilon)}{r_{\min}} I_{\overline p} + \Lambda^{-1}\right) = 1, 
  \end{align}
  and hence, 
  \begin{align}\label{eq:low-of-large-numbers-matrix}
    \mathrm{Pr}\left( \left( \frac{M}{r_{\min}}\Phi^\top \Phi  + \Lambda^{-1}\right)^{-1} \succ \left(\frac{NM(1+\epsilon)}{r_{\min}} I_{\overline p} + \Lambda^{-1}\right)^{-1}\right) = 1.
  \end{align}
  By combining \cref{eq:low-of-large-numbers-matrix} with \cref{eq:woodbury}, we can obtain 
  \begin{align}\label{eq:proof_thm1_bound1}
    \mathrm{Pr} \left( \widetilde\sigma^2(x) > \sum_{p \leq \overline{p}} \frac{r_{\min} \lambda_p}{r_{\min}+NM\lambda_p(1+\epsilon) }\ \phi_p(x)^2\right) = 1.    
  \end{align}
  Since $r_{\min} > 0$, 
  the right-hand side of the inequality inside the probability is bounded by
  \begin{align}
    \begin{aligned}
      \frac{r_{\min} \lambda_p}{r_{\min}+NM\lambda_p(1+\epsilon) } 
      &> \frac{r_{\min} \lambda_p}{r_{\min}+NM \lambda_p} (1-\epsilon).
    \end{aligned}\label{eq:proof_thm1_bound2}
  \end{align}
  By substituting \cref{eq:proof_thm1_bound2} into \cref{eq:proof_thm1_bound1}, we obtain
  \begin{align}
    \mathrm{Pr} \left( \widetilde\sigma^2(x) > \sum_{p \leq \overline{p}} \frac{r_{\min} \lambda_p}{r_{\min}+NM \lambda_p} (1-\epsilon)\ \phi_p(x)^2\right) = 1.
  \end{align}
  Integrating both sides of the expression inside the probability with measure $\nu(x)$ completes the proof.
  \qed

\subsection*{Proof of \cref{thm:err-bound}}
  First, by using the triangle inequality, we have
  \begin{align}
    |\overline r_{\min} - r_{\min}| \le |\overline r_{\min} - r(x_{i^*})| + |r(x_{i^*}) - r_{\min}|.
  \end{align}
  Similarly, when an event $\delta > \left|r\left(x_{i^*}\right)-r_{\min}\right|$ occurs,
  \begin{align}
    \begin{aligned}
      \left|\overline{r}_{\min}-r\left(x_{i^*}\right)\right| &\geq\left|\overline{r}_{\min}-r_{\min}\right|-\left|r\left(x_{i^*}\right)-r_{\min}\right|\\
      &>\left|\overline{r}_{\min}-r_{\min}\right|-\delta.
    \end{aligned}
  \end{align}
  Combining these inequalities, we obtain
  \begin{align}
    \begin{aligned}\label{eq:thm2_triangle} 
      & \operatorname{Pr}\left(\left|\overline{r}_{\min }-r_{\min }\right| \geq \epsilon, \delta>\left|r\left(x_{i^*}\right)-r_{\min}\right|\right) \\ & \quad \leq \operatorname{Pr}\left(\left|\overline{r}_{\min}-r_{\min}\right| \geq \epsilon,\left|\overline{r}_{\min}-r\left(x_{i^*}\right)\right| \geq\left|\overline{r}_{\min}-r_{\min}\right|-\delta\right) \\ & \quad \leq \operatorname{Pr}\left(\left|\overline{r}_{\min}-r\left(x_{i^*}\right)\right| \geq \epsilon-\delta\right).
    \end{aligned}
  \end{align}
  For the rightmost side, we have the following inequality: 
  \begin{align}\label{eq:thm2_almost_done}
    \begin{aligned}
      & \operatorname{Pr}\left(\left|\overline{r}_{\text {min }}-r\left(x_{i *}\right)\right| \geq \epsilon-\delta\right)\\
      & \quad\leq \operatorname{Pr}\left(\max _i\left|\overline{r}\left(x_i\right)-r\left(x_i\right)\right| \geq \epsilon-\delta\right) \\
      & \quad=\operatorname{Pr}\left(\exists i, \quad\left|\overline{r}\left(x_i\right)-r\left(x_i\right)\right| \geq \epsilon-\delta\right) \\
      & \quad=1-\operatorname{Pr}\left(\forall i \in\{1,2, \ldots, N\}, \quad\left|\overline{r}\left(x_i\right)-r\left(x_i\right)\right|<\epsilon-\delta\right) \\
      & \quad=1-\prod_{i=1}^N \operatorname{Pr}\left(\left|\overline{r}\left(x_i\right)-r\left(x_i\right)\right|<\epsilon-\delta\right).
      \end{aligned}          
  \end{align}
  We will find a bound for the probability on the rightmost side.
  For the Chebyshev's inequality~\cite{kallenberg1997foundations}
  \begin{align}
    \mathrm{Pr}(|Z|\ge \epsilon)\le \frac{\bbE[{|Z|^2}]}{\epsilon^2},
  \end{align}
  we substitute $Z=\overline r(x)-r(x)$, where $\overline r(x)$ denotes the unbiased sample variance of $\{u_j(x)\}_{j=1}^M$ and $r(x)$ denotes the variance of $u(x)$:
  \begin{align}
    \mathrm{Pr}(|\overline r(x)-r(x)|\ge \epsilon)\le \frac{\bbE[{|\overline r(x)-r(x)|^2}]}{\epsilon^2} = \frac{\bbV[\overline r(x)]}{\epsilon^2}.
  \end{align}
  The variance of the unbiased sample variance is known to be written as
  \begin{align}
    \bbV[\overline r(x)] = \frac{1}{M}\left(\mu_4[u(x)]-\frac{M-3}{M-1} r^2(x)\right),
  \end{align}
  where $\mu_4[u(x)]$ denote the 4th central moment of $u(x)$ (see \cite{Mood1950Introduction}).
  Note that $\mu_4[u(x)]=3r^2(x)$ when $u$ follows a normal distribution.
  By substituting this, we obtain the following concentration inequality:
  \begin{align}\label{eq:sample_variance_ineq}
    \mathrm{Pr}(|\overline r(x)-r(x)|\ge \epsilon) \le \frac{2r^2(x)}{\epsilon^2 (M-1)},
  \end{align}
  and thus, 
  \begin{align}
    \begin{aligned}
       \operatorname{Pr}(|\overline{r}(x)-r(x)|<\epsilon)
      &=1-\operatorname{Pr}(|\overline{r}(x)-r(x)| \geq \epsilon) \\
      &\geq \max \left\{1-\frac{2 r^2(x)}{\epsilon^2(M-1)}, 0\right\}.     
    \end{aligned}
  \end{align}
  Substituting this into \cref{eq:thm2_almost_done} yields
  \begin{align}
    \begin{aligned}
    \operatorname{Pr}\left(\left|\overline{r}_{\text {min }}-r\left(x_{i *}\right)\right| \geq \epsilon-\delta\right) 
    & \leq 1-\prod_{i=1}^N \max \left\{1-\frac{2 r^2\left(x_i\right)}{(\epsilon-\delta)^2(M-1)}, 0\right\}.
    \end{aligned}
  \end{align}
  Combining this and \cref{eq:thm2_triangle} completes the proof. \qed

\subsection*{Proof of \cref{col:conservative-bound}}
  Let $c \in(0, \infty)$ be a constant. 
  We have
  \begin{align}  
      & \forall \lambda_p \geq c, \quad \frac{\lambda_p}{r_{\min}+NM \lambda_p} \geq \frac{c}{r_{\min}+NM c}, \\
      & \forall \lambda_p<c, \quad \frac{\lambda_p}{r_{\min}+NM \lambda_p} \geq 0.
  \end{align}
  For $I_{\geq c}\coloneqq\left\{p \in I \mid \lambda_p \geq c\right\}$, we have
  \begin{align}
    L_{\mathrm{IMSE}} \geq r_{\min} \sum_{p \in I_{\geq c}} \frac{c}{r_{\min}+NM c}=\frac{c\ r_{\min}}{r_{\min}+NM c}\left|I_{\geq c}\right|.
  \end{align}
  Because $\lambda_p$ is independent of $N$ and $M$, substituting $N=1$ and $M=1$ into the above equation yields the boundedness of $\left|I_{\geq c}\right|$:
  \begin{align}
    \left|I_{\geq c}\right| \leq \frac{\left(r_{\min}+c\right)}{c\ r_{\min}} L_{\mathrm{IMSE}} <\infty.
  \end{align}
  Thus, there exists $C>0$ such that for all $N\in\bbN$ and $M\in\bbN$,
  \begin{align}
    L_{\mathrm{IMSE}} \geq \frac{c\ r_{\min}}{r_{\min}+NM c}\left|I_{\geq c}\right| \geq \frac{C}{NM},
  \end{align}
  which completes the proof.
  \qed

\subsection*{Review of Gaussian Process Regression}

We provide a brief review of the fundamentals of Gaussian process regression so that this paper is effectively self-contained, 
The main purpose is to give an overview of the eigenfunction expansions of kernels and Gaussian processes, which are necessary for the proof of \cref{thm:convergence-MSE-hetero-min}.

\begin{definition}[Positive-definite kernel]\label{def:PDK}
  For a set $\mathcal{X}\subseteq \bbR^d$, the function $k: \mathcal{X} \times \mathcal{X} \rightarrow \mathbb{R}$ is a positive-definite kernel if the following is satisfied for any $n \in \mathbb{N},\left[c_1, \ldots, c_n\right]^\top \in \mathbb{R}^n$ and $\left[x_1, \ldots, x_n\right]^\top \in \mathcal{X}^n$:
  \begin{align}
    \sum_{i=1}^n \sum_{j=1}^n c_i c_j k\left(x_i, x_j\right) \geq 0 .
  \end{align}
\end{definition}

\begin{definition}[Gaussian process]\label{def:GPR}
  For a set $\mathcal{X}\subseteq \bbR^d$, let $k: \mathcal{X} \times \mathcal{X} \rightarrow \mathbb{R}$ be a positive-definite kernel and $m: \mathcal{X} \rightarrow \mathbb{R}$ be a real-valued function.
  For a function $f: \mathcal{X} \rightarrow \mathbb{R}$, suppose that for any $n \in \mathbb{N}$ and $X=\left[x_1, \ldots, x_n\right]^\top \in \mathcal{X}^n$, the vector value $f_X=\left[f\left(x_1\right), \ldots, f\left(x_n\right)\right]^{\top}$ is a random variable that follows the multivariate Gaussian distribution $\mathcal{N}\left(m_X, k_{X X}\right)$ with mean $m_X=\left[m\left(x_1\right), \ldots, m\left(x_n\right)\right]^{\top}\in \mathbb{R}^{n}$ and covariance $k_{X X}=\left[k\left(x_i, x_j\right)\right]_{i, j=1}^n \in \mathbb{R}^{n \times n}$.
  Then, $f$ is called a Gaussian process $\mathcal{GP}(m, k)$ with mean $m$ and covariance $k$.
\end{definition}

\begin{definition}[reproducing kernel Hilbert space, RKHS]\label{def:RKHS}
  Consider a set $\mathcal{X}\subseteq \bbR^d$ and a positive-definite kernel $k$ on $\mathcal{X}$.
  A Hilbert space $\mathcal{H}_k$ consisting of functions on $\mathcal{X}$ whose inner product is defined by $\langle\cdot, \cdot\rangle_{\mathcal{H}_k}$ is called a reproducing kernel Hilbert space (RKHS) for a kernel $k$ if $\calH_k$ satisfies the following conditions:
  \begin{enumerate}
    \item For any $x \in \mathcal{X}$, $k(\cdot, x) \in \mathcal{H}_k$.
    \item For any $x \in \mathcal{X}$ and $f \in \mathcal{H}_k$, $f(x)=\langle f, k(\cdot, x)\rangle_{\mathcal{H}_k}$.
  \end{enumerate}
\end{definition}

Let $\nu$ be a Borel measure on $\mathcal{X}$.
Let $L_2(\nu)$ be a Hilbert space consisting of square-integrable functions on $\nu$.
Define the map $T_k: L_2(\nu) \rightarrow L_2(\nu)$ as follows:
\begin{align}\label{eq:t_k}
  T_k f\coloneqq\int k(\cdot, x) f(x) \rmd \nu(x), \quad f \in L_2(\nu) .
\end{align}
Since $T_k$ is a compact, positive-valued, self-adjoint function, eigenfunction expansion is possible according to spectral theory:
\begin{align}\label{eq:T_k-expansion}
  T_k f=\sum_{i \in I} \lambda_i\left\langle\phi_i, f\right\rangle_{L_2(\nu)} \phi_i,
\end{align}
where the convergence is in the $L_2(\nu)$-sense, and $I \subseteq \mathbb{N}$ denotes the index set.
For example, $I=\mathbb{N}$ if RKHS is infinite dimensional and $I=\{1, \ldots, K\}$ if $K(\in \mathbb{N})$ dimensional.
Let $\left\{(\phi_i, \lambda_i)\right\}_{i \in I} \subset L_2(\nu) \times(0, \infty)$ be the eigenfunction and eigenvalue of $T_k$, that is, 
\begin{align}\label{eq:eigen-equation}
  T_k \phi_i=\lambda_i \phi_i, \quad i \in I,\ \lambda_1 \geq \lambda_2 \geq \cdots>0 .
\end{align}
Note that the eigenfunctions $\left\{\phi_i\right\}_{i \in \mathbb{N}}$ satisfy
\begin{align}\label{eq:eigen-orthogonal}
  \int \phi_i(x) \phi_j(x) \rmd \nu(x) = \delta_{ij},
\end{align}
where $\delta_{i j}$ is the Kronecker delta.

Mercer's theorem states that the kernel $k$ can be represented by a system of eigenvalues and eigenfunctions.
\begin{proposition}[Mercer's theorem, \cite{Steinwart2008Support}]
  Let $\mathcal{X}\subseteq \bbR^d$ be a compact topological space, $k: \mathcal{X} \times \mathcal{X} \rightarrow \mathbb{R}$ be a continuous kernel, and $\nu$ be a bounded Borel measure on $\calX$.
  Then, 
  \begin{align}\label{eq:expansion}
    k\left(x, x^{\prime}\right)=\sum_{i \in I} \lambda_i \phi_i(x) \phi_i\left(x^{\prime}\right), \quad x, x^{\prime} \in \mathcal{X},
  \end{align}
  where the convergence is absolute and uniform over $x, x^{\prime} \in \mathcal{X}$.
\end{proposition}

Note that the eigenvalues and eigenfunctions $\left\{(\phi_i, \lambda_i)\right\}_{i \in I}$ in \cref{eq:expansion} depend on the measure $\nu$ because they are defined through integral operators of \cref{eq:t_k}.

\begin{definition}
  A kernel that has only a finite number of nonzero eigenvalues is called a degenerate kernel.
  A kernel that is not degenerate is called a nondegenerate kernel.
\end{definition}

\end{document}